\documentclass[11pt]{article}
\usepackage{graphicx}
\usepackage{amsmath,amsfonts,amssymb}
\usepackage{amsmath}
\usepackage{graphicx}
\usepackage{epsfig}
\usepackage{epstopdf}
\usepackage{color}
\def\Limsup{\mathop{{\rm Lim}\,{\rm sup}}}

\def\tto{\;{\lower 1pt \hbox{$\rightarrow$}}\kern -10pt
\hbox{\raise 2pt \hbox{$\rightarrow$}}\;}
\def\Hat{\widehat}

\def\ra{\rangle}
\def\la{\langle}
\def\ve{\varepsilon}
\def\B{\Bbb B}
\def\h{\hfill\Box}
\def\R{\Bbb R}
\def\N{\Bbb N}
\def\ox{\bar{x}}

\def\oy{\bar{y}}

\def\epi{\mbox{\rm epi}}

\def\dom{\mbox{\rm dom}}

\def\h{\hfill\square}
\def\dn{\downarrow}
\def\O{\Omega}
\def\ph{\varphi}

\def\vaep{\varepsilon}
\newcounter{lk}

\begin{document}
\begin{center}
{\bf Subgradients of Minimal Time Functions without Calmness}\\[1ex]
NGUYEN MAU NAM\footnote{Fariborz Maseeh Department of Mathematics and Statistics, Portland State University, PO Box 751, Portland, OR 97207, USA(mau.nam.nguyen@pdx.edu). The research of this author was partially supported by the NSF under grant DMS-1411817} and DANG VAN CUONG\footnote{Department of Mathematics, Duy Tan University, Da Nang, Vietnam (dvcuong@duytan.edu.vn).}
\end{center}
\small{\bf Abstract.}  In recent years there has been great interest in variational analysis of a class of nonsmooth functions called the minimal time function. In this paper we continue this line of research by providing new results on generalized differentiation of this class of functions, relaxing assumptions imposed on the functions and sets involved for the results. In particular, we focus on the singular subdifferential and the limiting subdifferential of this class of functions. \\[1ex]
{\bf Key words.} minimal time function, $\ve-$Fr\'echet subdifferential, limiting subdifferential, singular subdifferential, calmness.\\
\noindent {\bf AMS subject classifications.} 49J52, 49J53, 90C31

\newtheorem{Theorem}{Theorem}[section]
\newtheorem{Proposition}[Theorem]{Proposition}
\newtheorem{Remark}[Theorem]{Remark}
\newtheorem{Lemma}[Theorem]{Lemma}
\newtheorem{Corollary}[Theorem]{Corollary}
\newtheorem{Definition}[Theorem]{Definition}
\newtheorem{Example}[Theorem]{Example}
\renewcommand{\theequation}{\thesection.\arabic{equation}}
\normalsize
\section{Introduction}

It has been widely accepted that convex analysis is one of the most important
and useful areas of mathematical sciences, providing the mathematical foundation for convex optimization, a very fast growing field with many applications to many fields such as economics, computational statistics, compressed sensing and machine learning. Convex optimization builds effective numerical optimization algorithms to deal with both smooth and nonsmooth optimization problems involving large data sets encountered in many practical applications, especially in the recent time of \emph{big data}. At the same time, it is desirable to provide mathematical background and numerical optimizations for optimization problems in which the objective functions are both nonconvex and nonsmooth. This is the driving force for the development of nonsmooth/variational analysis. Started with the pioneering work of Clarke, Mordukhovich,  Rockafellar and others in the 1970's, variational analysis is now a mature area of mathematics; see \cite{bz,CL,CL1,m-book1,rw} and the references therein.

The class of distance functions is perhaps one of the most important examples of nonsmooth functions. Meanwhile, this class of functions plays a crucial role in many aspects of optimization. There has been extensive research on variational analysis of distance functions and their generalizations in the literature. In particular, the reader can find subdifferential formulas in the sense of convex analysis and the Clarke nonconvex subdifferential in \cite{BFQ}, while  the Fr\'echet subdifferetial formula was obtained in \cite{BT} and the limiting subdifferential formula was obtained in \cite{bmn05}.  One of the most natural generalizations of the distance function is the minimal time function, obtained by replacing the norm function that defines the distance function by a Minkowski gauge. Recall that given a nonempty closed set $\O$ in a normed space $X$ and a non empty closed bounded convex set $F$, the minimal time function to $\O$ with the \emph{constant dynamics} $F$ is given by
\begin{equation}\label{mt}
T_{\O}^F(x):=\inf\big\{ t\geq 0\; \big|\; (x+tF)\cap \O\neq\emptyset\big\},\; x\in X.
\end{equation}
It turns out that the minimal time function has the representation
\begin{equation*}
T_{\O}^F(x):=\inf\big\{ \rho_F(w-x)\; \big|\; w\in \O\big\},\; x\in X,
\end{equation*}
where $\rho_F(u):=\inf\{t\geq 0\; |\; u\in tF\}$ for $u\in X$.  Generalized differentiation in both convex and nonconvex setting for the class of minimal time functions was considered in \cite{CW04a,CW04b}. Further study in Banach spaces was presented in several research papers including \cite{HNg,JH,Meng,singbmn,SY}.

This paper concerns new results of variational analysis of the minimal time function. In particular, we focus on obtaining the singular and the limiting subdifferential formula for this class of functions. The result obtained in our paper extend the line of research in this direction by providing various subdifferential formulas for the minimal time function without requiring the calmnes  as initiated in \cite{SY}.

Throughout this paper we consider a real normed space $X$ with a given norm $\|\cdot\|$. The dual space of $X$ is denoted by $X^*$ and the paring of an element $x^*\in X^*$ and $x\in X$ is denoted by $\la x^*,x\ra$, i.e., $\la x^*, x\ra:=x^*(x)$.  We always assume that $F$ is nonempty closed bounded convex set and $F\ne \{0\}.$  The closed ball centered at $\ox$ with radius $r>0$ is denoted by $\B(\ox; r),$  the open ball centered at $\ox$ with radius $r>0$ is denoted by $\B^{\rm o}(\ox; r),$ and the closed unit balls of $X$ and $X^*$ are denoted by $\B$ and $\B^*$, respectively.

\section{Preliminaries}
In this section we present basic notions and results of variational analysis used throughout the paper. The readers are referred to the books \cite{bz,CL,CL1,m-book1} for more details.

Given an extended real-valued function $f\colon X\to (-\infty, \infty]$, with the domain $\dom(f):=\{x\in X\; |\; f(x)<\infty\}$, and given $\ve\geq 0$, the $\ve-$\emph{Fr\'echet subdifferential} (or the set of $\ve-$Fr\'echet subgradients) of $f$ at a point $\ox\in \dom(f)$ is defined by
\begin{equation*}
\Hat\partial_\ve f(\ox):=\big\{x^* \in X^*\;\big |\; \liminf_{x\to \ox}\dfrac{f(x)-f(\ox)-\la x^*, x-\ox\ra}{\|x-\ox\|}\geq -\vaep\big\}.
 \end{equation*}
If $\ox\notin \dom(f)$, we set $\Hat\partial_\ve f(\ox)=\emptyset$. In the case where $\ve=0$, we use the notation $\Hat\partial f(\ox)$ instead of $\Hat\partial_0f(\ox)$ for simplicity. If $f$ is a convex function, the $\ve-$Fr\'echet subdifferential has a simpler representation, namely,
\begin{equation*}
\Hat\partial_\ve f(\ox):=\big\{x^*\in X^*\;\big |\; \la x^*, x-\ox\ra \leq f(x)-f(\ox)+\vaep \|x-\ox\|\; \mbox{\rm for all } x\in X\big\},
\end{equation*}
which reduces to the classical subdifferential in the sense of convex analysis when $\ve=0$.

Based on the $\ve-$Fr\'echet subdifferential, two major concepts of variational analysis called the \emph{singular subdifferential} and the \emph{limiting subdifferential} are defined using the \textit{sequential Painlevé-Kuratowski upper limit} as follows:
\begin{align}\label{Mordukhovichsubdifferential}
\partial f(\bar{ x}):=\Limsup_{x\xrightarrow{f}\ox}\Hat\partial_\vaep f(x).
\end{align}
and
\begin{equation}\label{singsub}\partial^{\infty}f(\ox)=\Limsup_{x \xrightarrow{f} \bar x,\varepsilon,\lambda \downarrow 0}\lambda \Hat{\partial}_{\vaep}f(x).\end{equation}
Here $x \xrightarrow{f} \bar x$ means that $x\rightarrow \ox$ and $f(x)\rightarrow f(\ox)$.

The inclusion $\widehat\partial f(\bar x) \subset \partial f(\bar x)$ is valid for any $\bar x \in X.$ If $f$ is convex, then
\begin{equation*}\label{subdifferential_convex_case}\widehat\partial f(\bar x) =\partial f(\bar x)=\{x^* \in X^* \mid \langle x^*, x- \bar x \rangle \le f(x)-f(\bar x) \; \mbox{\rm for all }  x \in X\},\end{equation*}
i.e., the Fr\'{e}chet subdifferential and the Mordukhovich subdifferential of $f$ at $\bar x$ coincide with the subdifferential of $f$ at $\bar x$ in the sense of convex analysis.

Both subdifferential notions \eqref{Mordukhovichsubdifferential} and \eqref{singsub} have geometric representations in terms of normal cones to sets defined in what follows. Given a subset $\Omega \subset X$, we use the notation $x \xrightarrow{\Omega}u$ to mean that $x \rightarrow u$ and $x\in \Omega$. For any $x \in \Omega$ and $\varepsilon \geq 0$, the set of \textit{ $\varepsilon$-normals} to $\Omega$ at $x$ is defined by
\begin{equation*}\label{varepsilon-normals}
\widehat N_\varepsilon(x;\Omega):=\big\{ x^*\in X^*\; \big |\;  \limsup_{u \xrightarrow{\Omega}x} \dfrac{\langle x^*, u-x \rangle}{\|u-x\|} \leq \varepsilon \big\}.
\end{equation*}
The set $\widehat N(x;\Omega):=\widehat N_0(x;\Omega)$ is called the \textit{Fr\'echet normal cone} to $\Omega$ at $x.$ If $x \not\in \Omega$, we put $\widehat N_\varepsilon(x;\Omega):=\emptyset$ for all $\varepsilon \geq 0.$

Given $\bar x \in \Omega$, the set
\begin{equation*}\label{1.2}
N(\bar x; \Omega) :=\Limsup_{x \rightarrow \bar x,\varepsilon \downarrow 0}\,\widehat N_\varepsilon(x; \Omega)
\end{equation*}
is called the \textit{Mordukhovich normal cone} or the \textit{limiting  normal cone} to $\Omega$ at $\bar x$.  We put $N(\bar x;\Omega)=\emptyset$ if $\bar x \not\in \Omega.$

It is clear that $\widehat N(x; \Omega) \subset N(x ; \Omega)$ for all $x \in \Omega$. In the case where  $\Omega$ is a convex, one has the following simple representation:
$$ \widehat N_\varepsilon(\bar x;\Omega)=\big\{ x^*\in X^*\; \big |\;  \langle x^*, x- \bar x \rangle  \leq \varepsilon \|x-\bar x\|\; \mbox{\rm for all } x \in \Omega\big\}$$
for all $\varepsilon \geq 0$ and $\bar x \in \Omega.$ Moreover, both $\widehat N(\bar x; \Omega)$ and $N(\ox;\O)$ coincide with the convex cone to $\Omega$ at $\bar x$ in the sense of convex analysis, that is,
\begin{align*}\label{normal_cone_convex_analysis}
\widehat N(\bar x; \Omega)=N(\ox;\O)=\big\{x^*\in X^* \; \big |\;  \langle x^*, x-\bar x \rangle \leq 0  \; \mbox{\rm for all } x \in \Omega\big\}.\end{align*}

%
%

\section{Fr\'echet Singular Subgradients}

In this section we introduce and study the Fr\'echet singular subdifferential of extended real-valued functions. In addition, we present a new result on the singular subdifferential in the convex case in Banach spaces, while a similar result for the limiting singular subdifferential is well-known in Aspund spaces; see \cite{m-book1}.

\begin{Definition} Let $f\colon X\to (-\infty, \infty]$ be an extended real-valued function and let $\ox\in \dom(f)$. Define the Fr\'echet singular subdifferential of $f$ at $\ox$ by
\begin{equation*}
\Hat\partial^\infty f(\ox):=\big\{x^*\in X^*\;\big |\; (x^*, 0)\in \Hat N((\ox, f(\ox)); \epi(f))\big\}.
\end{equation*}
\end{Definition}
In the theorem below, we study this concept in connection with the limiting singular subdifferential and the convex normal cone to the domain of the function $f$ involved.

The proof of the proposition below is straightforward.

\begin{Proposition}\label{sn}
Let $f:X\to(-\infty,\infty]$  be a convex function and let $\ox \in \dom(f),$ where $X$ is  a normal space. Then
\begin{equation*}
\Hat\partial^\infty f(\ox)=\Hat N(\ox; \dom(f))=N(\ox; \dom(f)).
\end{equation*}
\end{Proposition}

Now we are ready to prove the main result of this section concerning the Fr\'echet singular subdifferential in the convex case.
\begin{Theorem} Let $X$ be a Banach space and let $f\colon X\to (-\infty, \infty]$ be a l.s.c convex function. Then
\begin{equation*}
\partial^\infty f(\ox)=\Hat\partial^\infty f(\ox)=\Hat N(\ox; \dom(f))=N(\ox; \dom(f)).
\end{equation*}
In addition, $x^*\in \Hat\partial^\infty f(\ox)$ if and only if there exist $x_k\to \ox$, $f(x_k)\to f(\ox)$, $\lambda_k\dn 0$, $x_k^*\in \partial f(x_k)$ and $\lambda_k x^*_k\xrightarrow{\|\cdot\|} x^*$.
\end{Theorem}
{\bf Proof.} By Proposition \ref{sn}, we have
\begin{equation*}
\Hat\partial^\infty f(\ox)=\Hat N(\ox; \dom(f))=N(\ox; \dom(f)).
\end{equation*}
Let us now prove that $N(\ox; \dom(f))\subset \partial^\infty f(\ox)$. Fix any $x^* \in N(\ox; \dom(f))$. Choose a sequence $\ve_k\dn 0$. By \cite[Proposition~3.15]{ph}, there exists $u^*_k\in X^*$ such that
\begin{equation*}
\la u^*_k, x-\ox\ra \leq f(x)-f(\ox)+\ve_k\; \mbox{\rm for all }x\in X.
\end{equation*}
We have
\begin{equation*}
\langle {x^*},x - {\rm{ }}\bar x\rangle  \le 0\; \mbox{\rm for all } x \in \dom(f),
\end{equation*}
and hence  for such $x$,
\begin{equation*}
\langle k(\left\| {u_k^*} \right\| + 1){x^*},x - {\rm{ }}\bar x\rangle  \le 0.
\end{equation*}
This implies
\[\langle u_k^* + k(\left\| {u_k^*} \right\| + 1){x^*},x - {\rm{ }}\bar x\rangle  \le f(x) - f(\bar x)+ {\varepsilon _k}\;\mbox{\rm for all}\; x \in X.\]
Applying \cite[Theorem~3.1.1]{z} with $\beta=1$, for every $k\in \N$ there exist $x_k$, $e^*_k\in \B^*$,  and $|\lambda_k|\leq 1$ with
\begin{align*}
&\| {x_k - \bar x} \| \le  \sqrt{\varepsilon_k},\\
&u_k^* + k(\left\| {u_k^*} \right\| + 1){x^*}+\sqrt{\ve_k}\,(e^*_k+\lambda_k(u_k^* + k(\left\| {u_k^*} \right\| + 1){x^*}))\in\partial f(x_k),\\
&|f(x_k)-f(\ox)|\leq \sqrt{\ve_k}+\ve_k.
\end{align*}
Let $z^*_k:=e^*_k+\lambda_k(u_k^* + k(\left\| {u_k^*} \right\| + 1){x^*})$. So there exist $\tilde{x}_k^*\in \partial f(x_k)$  such that
\begin{equation*}
u_k^* + k(\left\| {u_k^*} \right\| + 1){x^*}+ \sqrt {{\varepsilon_k}} z_k^* = \tilde{x}_k^*.
\end{equation*}
Then
\[\frac{\tilde{x}_k^*}{{k(\left\| {u^*_k} \right\| + 1)}} = {x^*} + \frac{{u^*_k}}{{k(\left\| {u^*_k} \right\| + 1)}} +\frac{{\sqrt {{\varepsilon _k}}\, z_k^*}}{{k(\left\| {u^*_k} \right\| + 1)}} \to x^*\; \mbox{\rm as}\;k \to \infty .\]
Note that $\tilde{x}_k^* \in \partial f(x_k)$ and ${\lambda' _k}:= \frac{1}{{k(\left\| {u_k^*} \right\| + 1)}} \to {0^ + }\; \mbox{\rm as}\;k \to \infty$. In addition, $f(x_k)\to f(\ox)$ as $k\to \infty$, so $x^*\in \partial^\infty f(\ox)$.

It remains  to show that $\partial ^\infty f(\ox)\subset N(\ox; \dom(f))$. Fix any $x^* \in \partial ^\infty f(\ox)$. Then there exist ${\lambda _k} \to {0^ + }, {x_k} \to \ox, f(x_k)\to f(\ox), x_k^* \in \partial f({x_k})$ and ${\lambda _k}{x_k^*}\mathop  \to \limits^ {w^*}  {x^*}$.\\
Since $x_k^* \in \partial f({x_k})$, we have
\begin{equation*}
\langle x_k^*,x - {x_k}\rangle  \le f(x) - f({x_k})\; \mbox{\rm for all}\;  x \in X.
\end{equation*}

Hence,
\begin{equation*}
\langle {\lambda _k}x_k^*,x - {x_k}\rangle  \le {\lambda _k}(f(x) - f(\bar x)\; \mbox{\rm for all } x \in X.
\end{equation*}
Letting $k \to \infty $, we have $\langle {x^*},x - \bar x\rangle  \le 0$ for all $x\in X$, and  thus ${x^*} \in N(\bar x; \dom(f))$, which completes the proof.$\h$

\section{Fr\'echet Singular Subgradients of  Minimal Time Functions}

In this section we study Fr\'echet singular subdifferential formulas for the minimal time function in both in-set and out-of-set settings.

Following \cite{singbmn}, we define the following sets:
$$S_{\vaep}^*:=\{x^*\in X\; |\;\ 1-\vaep\|F\|\leq\sigma_{F}(-x^*)\leq 1+\vaep\|F\|\} \;\mbox{\rm for } \ve\geq 0,\ S^*:=S_0^*,$$
$$C^*:=\{x^*\in X\; |\; \sigma_F(-x^*)\leq 1\},\ F^*_+:=\{x^*\in X^*\; |\; \la x^*,q\ra\geq 0\ \text{for all}\ q\in F\},$$
where
$$\sigma_F(x^*):=\sup_{x\in F}\la x^*,x\ra\;\text{for}\;  x^*\in X^*,\ \|F\|:=\sup\{\|q\|\; |\;  q\in F\}.$$

Given $r>0$, define the enlargement set
\begin{equation*}
\O_r:=\big\{x\in X\; \big |\; T_{\O}^F(x)\leq r\big\}.
\end{equation*}

In the next proposition we present some basic facts about the minimal time function \eqref{mt}. The reader can find the detailed proof in \cite{JH,singbmn}.

\begin{Proposition} Consider the minimal time function \eqref{mt}. The following properties hold:\\[1ex]
{\rm\bf (i)} For $x\in X$, $T_{\O}^F(x)=0$ if and only if $x\in \O$. \\
{\rm\bf (ii)} For any $x \in \O_r$ with $r>0$ and $t\geq 0$ we have
\begin{equation*}
T_{\O}^F(x-tq)\leq r+t\; \mbox{\rm whenever }q\in F.
\end{equation*}
{\rm\bf (iii)} If $x\notin \O_r$ with $T_{\O}^F(x)<\infty$, then
\begin{equation*}
T_{\O}^F(x)=T_{\O_r}^F(x)+r.
\end{equation*}
\end{Proposition}

Let us now present a formula for computing the Fr\'echet singular subdifferential of the minimal time function when the reference point is in the target set.

\begin{Proposition}\label{p1} For any $\ox\in \O$, we have
\begin{equation*}
\Hat\partial^\infty T_{\O}^F(\ox)=\Hat N(\ox; \O)\cap F_{+}^*.
\end{equation*}
\end{Proposition}
{\bf Proof.} Fix any $x^*\in \Hat\partial^\infty T_{\O}^F(\ox)$ and let $\oy:=T_\O^F(\ox)=0$. Then $(x^*, 0)\in \Hat N((\ox, \oy); \epi(T_\O^F))$. Given any $\ve>0$, there exists $\delta>0$ such that
\begin{equation*}
\la x^*, x-\ox\ra \leq \ve(\|x-\ox\|+\lambda)\; \mbox{\rm whenever }\|x-\ox\|<\delta,  T_\O^F(x)\leq \lambda <\delta.
\end{equation*}
For any $x\in \O$ with $\|x-\ox\|<\delta$, let $\lambda=T_\O^F(x)=0$ and get
\begin{equation*}
\la x^*, x-\ox\ra\leq \ve \|x-\ox\|.
\end{equation*}
This implies $x^*\in \Hat N(\ox; \O)$. Now fix any $q\in F$ and $t>0$ sufficiently small such that $\ox-tq\in \B(\ox; \delta)$. Then $T_\O^F(\ox-tq)\leq t$, and so
\begin{equation*}
\la x^*, (\ox-tq)-\ox\ra\leq \ve(\|(\ox-tq)-\ox\|+t).
\end{equation*}
This implies $-\la x^*,q\ra\leq \ve(\|q\|+1)$. Letting $\ve\dn 0$ gives $\la x^*, q\ra\geq 0$, and so $x^*\in F_+^{*}$.

Let us prove the opposite inclusion. Fix any $x^*\in \Hat N(\ox; \O)\cap F_+^{*}.$ For any $\ve>0$, we can choose $\delta>0$ such that
\begin{equation*}
\la x^*, x-\ox\ra\leq \ve\|x-\ox\|\; \mbox{\rm whenever }\|x-\ox\|<\delta, x\in \O.
\end{equation*}
Fix any $x\in X$ with $\|x-\ox\|<\delta/2$ and fix any $t>0$ and $T_\O^F(x)\leq t<\delta$. Then we can find $0\leq t^\prime$ with $t^\prime\|F\|<t$ and $q\in F$ with $x+t^\prime q\in \O\cap \B(\ox; \delta)$. Thus,
\begin{equation*}
\la x^*, x+tq-\ox\ra\leq \ve(\|x+tq-\ox\|)\leq \ve(\|x-\ox\|+t^\prime\|q\|).
\end{equation*}
This implies
\begin{equation*}
\la x^*, x-\ox\ra \leq \ve(\|x-\ox\|+t).
\end{equation*}
Therefore, $x^*\in \Hat\partial^\infty T_\O^F(\ox)$. $\h$

Let us now consider the case where $\ox\notin\O$.

\begin{Theorem} Consider $\ox\notin\O$ and let $r:=T_\O^F(\ox)$. We have
\begin{equation*}
\Hat\partial^\infty T_\O^F(\ox)=\Hat N(\ox; \O_r)\cap F_{+}^*.
\end{equation*}
\end{Theorem}
{\bf Proof.} Fix any $x^*\in \Hat\partial^\infty T_\O^F(\ox)$. For any $\vaep>0$, find $\delta>0$ such that
\begin{equation*}
\la x^*, x-\ox\ra\leq \ve (\|x-\ox\|+|\lambda-r|)
\end{equation*}
whenever $\|x-\ox\|<\delta$ and $|\lambda -r|<\delta$, $\lambda\geq r$. Now we fix any $x\in \O_r$ with $\|x-\ox\|<\delta$. Then $T^F_\O(x)\leq r$, and so $(x, r)\in \epi(T^F_\O)$. Applying the inequality above with $\lambda:=r=T^F_\O(\ox)$ yields
\begin{equation*}
\la x^*, x-\ox\ra \leq \ve \|x-\ox\|.
\end{equation*}
Thus $x^*\in \Hat N(\ox; \O)$. Let us now show that $x^*\in F_+^*$. Fix any $q\in F$ and choose $t>0$ sufficiently small such that $\ox-tq\in \B(\ox; \delta)$. Then $T^F_\O(\ox-tq)\leq r +t$, and so $(\ox-tq, r+t)\in \epi (T^F_\O)$. With sufficiently small $t>0$, one has
\begin{equation*}
\la x^*, (\ox-tq)-\ox\ra\leq \ve (\|(\ox-tq)-\ox\|+|(r+t)-r|).
\end{equation*}
This implies $\la x^*, q\ra\geq 0,$ and hence $x^*\in F_+^*.$

Let us now proof the opposite inclusion. Fix any $x^*\in \Hat N(\ox; \O_r)\cap F_{+}^*.$ Then by Proposition \ref{p1}, $x^*\in \Hat\partial^\infty T^F_{\O_r}$.  In addition, $x^*\in \Hat N(\ox; \O_r)$, so for any $\ve>0$, there exists $\delta>0$ such that
\begin{equation*}
\la x^*, x-\ox\ra\leq \ve\|x-\ox\|\; \mbox{\rm whenever }\|x-\ox\|<\delta, x\in \O_r
\end{equation*}
and
\begin{equation*}
\la x^*, x-\ox\ra \leq \ve (\|x-\ox\|+\gamma)\; \mbox{\rm whenever }|x-\ox|<\delta, 0\leq \gamma<\delta, \gamma>T_{\O_r}^F(x).
\end{equation*}

Let us now fix any $x\in X$ and $\lambda\in \R$ with $\|x-\ox\|<\delta$ and $\lambda <\delta$ with $\lambda\geq T^F_\O(x)$. Consider the first case where $T^F_\O(x)\leq r$. Then $x\in \O_r$, and so
\begin{equation*}
\la x^*, x-\ox\ra\leq \ve\|x-\ox\|\leq \ve (\|x-\ox\|+|\lambda-T^F_\O(\ox)|).
\end{equation*}
Consider the second case where $\lambda\geq T^F_\O(x)>r=T_{\O}^F(\ox)$. Then $\lambda\geq T^F_\O(x)=r+T^F_{\O_r}(x)$, and so $T_{\O_r}^F(x)\leq \lambda-r<\delta$. Thus,
\begin{equation*}
\la x^*, x-\ox\ra \leq \ve (\|x-\ox\|+\lambda-r)=\ve (\|x-\ox\|+|\lambda-T^F_\O(\ox)|).
\end{equation*}
Therefore, $x^*\in \Hat\partial^\infty T^F_\O(\ox)$. $\h$

\section{$\ve-$Fr\'echet and Limiting Subgradients of  Minimal Time Functions}

In this section we study $\ve-$Fr\'echet and limiting subdifferentials of the minimal time function without imposing the calmness condition. We focus on the case where the reference point is outside of the target set as the other case has been considered in \cite{singbmn}.

\begin{Lemma}\label{Lmminimaltime} Consider the function
$$f(t):=\frac{at+a+c}{-(1+ab)t+1-2bc},$$
where $a,b,$ and $c$ are positive real numbers satisfying $c<\frac{1}{2b}$ and $1-2ab>0.$ Then $f$ is increasing on the interval $(0, \alpha),$ where $\alpha:=\frac{1-2ab}{1+ab}.$
\end{Lemma}
{\bf Proof.} Obviously,
$$f^\prime(t)=\frac{2a-abc+a^2b+c}{[-(1+ab)t+1-2bc]^2}=\frac{a(2-bc)+a^2b+c}{[-(1+ab)t+1-2bc]^2}>0$$
for all $0<t<\frac{1-2ab}{1+ab},$ and hence the lemma has been proved. $\h$

Now, we study  the $\vaep$-Fr${\rm \acute{e}}$chet subdifferential of minimal time function at points outside $\O$. The theorem below improves a result in \cite{singbmn} by removing the calmness assumption. We follow the proof from \cite{singbmn,SY}.
\begin{Theorem}\label{Theoepsub} Let $\ox\notin\O$ and $r:=T_{\O}^F(\ox)<\infty.$ Then for any $x^*\in \Hat N_{\vaep}(\ox;\O_r)\cap S^*_{\vaep}$ and $\vaep\geq0$ satisfying $1-2\vaep\|F\|>0,$ there exists a constant $\ell:=1+2\kappa\|F\|$ with $\kappa>\|x^*\|$ such that $x^*\in\Hat{\partial}_{\ell\vaep}T_{\O}^F(\ox).$
\end{Theorem}
{\bf Proof.}
Let $x^*\in \Hat N_{\vaep}(\ox;\O_r)\cap S_{\vaep}^*.$ Fix $k_0$ satisfying
 $$0<k_0<\frac{1-2\vaep\|F\|}{1+\|x^*\|\|F\|},$$
 and  set
$$\kappa:=\frac{\|x^*\|(k_0+1)+\vaep}{1-2\vaep \|F\|-k_0(1+\|x^*\|\|F\|)}.$$
It is easy to see that $\kappa>\|x^*\|$ is a constant. We will show that
\begin{equation}\label{eq1Theo2}
 \liminf_{x\to \ox}\frac{T_{\O}^F(x)-T_{\O}^F(\ox)-\la x^*, x-\ox\ra}{\|x-\ox\|}\geq -\ell \vaep,
\end{equation}
where $\ell:=1+2\kappa\|F\|.$ Using the proof of \cite[Proposition~4.6]{singbmn}, we only need to consider the case where $T_{\O}^{F}(x)=q< r,$ where $\|F\|>0.$ In addition, it suffices to consider the case where $\ve>0$ because the other case has been  considered in \cite{SY}. For any $\eta>0,$ suppose that  $\eta<2\vaep \|F\|.$ Let $\eta_0>0$ such that
\begin{equation}\label{eq1Theoeptime}0<\eta_0<\min\left\{k_0,\frac{\eta}{\|x^*\|}\right\}.
\end{equation}
Since $0<\eta_0<k_0,$ it follows from Lemma \ref{Lmminimaltime} with $a=\|x^*\|, b=\|F\|$ and $c=\vaep$ that
\begin{equation}\label{eq2Theo2}
0<\frac{\|x^*\|(\eta_0+1)+\vaep}{1-2\vaep \|F\|-\eta_0(1+\|x^*\|\|F\|)}<\frac{\|x^*\|(k_0+1)+\vaep}{1-2\vaep \|F\|-k_0(1+\|x^*\|\|F\|)}=\kappa.
\end{equation}

We first show that
\begin{equation}\label{eq2'Theo2}
  T_{\O}^F(\ox)-T_{\O}^F(x) \leq \kappa\|x-\ox\|.
\end{equation}
Since $x^*\in \Hat N_{\vaep}(\ox;\O_r)$ and $x^*\ne0,$ there exists $\delta>0$ such that
\begin{equation}\label{eq3Theoeptime}
\la x^*,x-\ox\ra<(\vaep +\eta_0\|x^*\|)\|x-\ox\|\leq (\vaep+\eta)\|x-\ox\|\ \text{for every} \ x\in\O_r\cap\mathbb B^{\rm o}(\ox,\delta).
\end{equation}

Let $\delta_1:=\frac{\delta}{2(1+\kappa\|F\|)}.$  Fix any $x\in\mathbb B^{\rm o}(\ox,\delta_1)$ such that $q:=T_{\O}^F(x)<r.$ Since $\delta_1<\delta,$ we have $x\in \mathbb B^{\rm o}(\ox,\delta). $

Since $\sigma_F(-x^*)\geq1-\vaep\|F\|,$  there exists $f\in F$ such that
\begin{equation}\label{eq7Theoeptime}  \la -x^*,f\ra >1-\vaep\|F\|-\eta_0.
\end{equation}
Take $z_t:=x-tf$ for $t>0.$ We claim that there exists $\hat t$ such that
\begin{equation}\label{clamihatt} T_{\O}^F(z_{\hat t})\geq r\ \text{and}\ 0<\hat t<\kappa\|x-\ox\|.
\end{equation}
Indeed, if $0<t<\min\{\frac{\delta}{2\|F\|},r-q\},$ then Proposition 4.1  implies that
$$z_t\in\O_r\cap\mathbb B^{\rm o}(\ox,\delta),$$
and thus we obtain from (\ref{eq3Theoeptime}) that
\begin{equation}\label{eq8Theoeptime}\la x^*,z_t-\ox\ra < (\vaep+\eta_0\|x^*\|)\|z_t-\ox\|.
\end{equation}
On the other hand, by (\ref{eq7Theoeptime}) and (\ref{eq1Theoeptime}) we have
$$\begin{aligned}\lim_{t\to+\infty}\frac{\la x^*,z_t-\ox\ra}{t}&=\la -x^*,f\ra>1-\vaep\|F\|-\eta_0\\
&>(\eta_0\|x^*\|+\vaep)\|F\|\geq (\eta_0\|x^*\|+\vaep)\|f\|=\lim_{t\to+\infty}\frac{(\eta_0\|x^*\|+\vaep)\|z_t-\ox\|}{t}.
\end{aligned}$$
So, we obtain for large enough $t$,
\begin{equation}\label{eq9Theoeptime}\la x^*,z_t-\ox\ra >(\vaep+\eta_0\|x^*\|)\|z_t-\ox\|.
\end{equation}
Combining (\ref{eq9Theoeptime}) with (\ref{eq8Theoeptime}) yields that there exists $\hat t>0$ such that
\begin{equation}\label{eq10Theoeptime} \la x^*,z_{\hat t}-\ox\ra=(\eta_0\|x^*\|+\vaep)\|z_{\hat t}-\ox\|.
\end{equation}
It follows from (\ref{eq8Theoeptime}) that
\begin{equation}\label{eq11Theoeptime} z_{\hat t}\notin \O_r\cap \mathbb B^{\rm o}(\ox,\delta) \ \text{or}\ z_{\hat t}=\ox\ \text{or}\  z_{\hat t}=\ox.
\end{equation}

By (\ref{eq7Theoeptime}) and (\ref{eq10Theoeptime}), we have
$$
\begin{aligned}
-\|x^*\|\|x-\ox\|+\hat t(1-\vaep\|F\|-\eta_0)&\leq \la x^*,x-\ox\ra +\hat t(1-\vaep\|F\|-\eta_0)\\
&\leq \la x^*,x-\ox\ra +\hat t\la x^*,-f\ra\\
&=\la x^*,z_{\hat t}-\ox\ra =(\eta_0\|x^*\|+\vaep)\|z_{\hat t}-\ox\|\\
&\leq (\eta_0\|x^*\|+\vaep)(\|x-\ox\|+\hat t \|f\|)\\
&\leq (\eta_0\|x^*\|+\vaep)\|x-\ox\|+(\eta_0\|x^*\|+\vaep)\hat t\|F\|,
\end{aligned}
$$
and, from (\ref{eq2Theo2}),  hence
\begin{equation}\label{eq12Theoeptime}0<\hat t<\frac{\|x^*\|(\eta_0+1)+\vaep}{1-2\vaep \|F\|-\eta_0(1+\|x^*\|\|F\|)} \|x-\ox\|<\kappa\|x-\ox\|.
\end{equation}
Since $x\in \mathbb B^{\rm o}(\ox,\delta_1),$ that is, $\|x-\ox\|<\frac{\delta}{2(\|F\|\kappa+1)},$ then
$$\|z_{\hat t}-\ox\|\leq \|x-\ox\|+\hat t\|F\|\leq (\kappa\|F\|+1)\|x-\ox\|<\frac{\delta}{2}<\delta.$$
It follows from (\ref{eq11Theoeptime}) that
$$z_{\hat t}\notin\O_r \ \text{or}\ z_{\hat t}=\ox $$
which implies that $T_{\O}^F(z_{\hat t})\geq r.$ Therefore, by using  (\ref{eq12Theoeptime}), there exists $\hat t$ satisfying (\ref{clamihatt}).

By Proposition 4.1 (ii), we have
$$T_{\O}^F(z_{\hat t})-T_{\O}^F(x)=T_{\O}^F(x-\hat tf)-T_{\O}^F(x)\leq \hat t.$$
It follows from (\ref{clamihatt}) that
$$
T_{\O}^F(\ox)-T_{\O}^F(x)=r-T_{\O}^F(x)\leq T_{\O}^F(z_{\hat t})-T_{\O}^F(x)\leq \hat t\leq  \kappa\|x-\ox\|,$$
and (\ref{eq2'Theo2}) was proved.

To obtain (\ref{eq1Theo2}), we only need to repeat the proof of \cite[Theorem 4.6]{singbmn}. We include the details for the convenience of the reader.
Since $\sigma_{F}(-x^*)\geq1-\vaep\|F\|,$ there exists $f\in F$ such that
\begin{equation}\label{eq15Theoeptime} \la -x^*,f\ra > 1-\vaep\|F\| -\eta.
\end{equation}

Take a sequence $\nu_k\downarrow0$ as $k\to\infty$. For any $k\in \mathbb N$, find $t_k\geq 0,$ $w_k\in \O,$ and $f_k\in F$ satisfying
$$q\leq t_k\leq q+\nu_k\ \text{and}\ w_k=x+t_kf_k.$$
Observe that
$$w_k=x-(r-t_k)f+(r-t_k)f+t_kf_k\subset x-(r-t_k)f+rF$$
when $k$ is sufficiently large. Thus for such $k$ we have
$$T_{\O}^F(x_k)\leq r\ \text{with}\ x_k:=x-(r-t_k)f.$$
Using (\ref{eq2'Theo2}) and the definition of $\delta_1$, we arrive subsequently at the upper estimates
\begin{equation}\label{eq4theo2}
  \begin{aligned} \|x_k-\ox\|&\leq \|x-\ox\|+(r-t_k)\|f\|\leq \|x-\ox\|+(r-q)\|F\|\\
  &\leq \|x-\ox\|+\kappa\|x-\ox\|\cdot\|F\|\leq (1+\kappa\|F\|)\delta_1<\delta,
  \end{aligned}
\end{equation}
and thus $x_k\in\B(\ox,\delta)$ for all $k$ is sufficiently large. Plugging now $x:=x_k$ into (\ref{eq3Theoeptime}) and employing the meddle estimate in (\ref{eq4theo2}), we get
$$
\begin{aligned}\la x^*,x-\ox\ra -(r-t_k)\la x^*,f\ra &\leq (\vaep +\eta)\|x_k-\ox\|\\
&\leq (\vaep+\eta)(1+\kappa\|F\|)\|x-\ox\|
\end{aligned}
$$
for the point $x$ fixed above. Letting $k\to\infty$ and using (\ref{eq15Theoeptime}), one has
$$
\begin{aligned}\la x^*,x-\ox\ra &\leq (r-q)\la x^*,f\ra +(\vaep+\eta)(1+\kappa\|F\|)\|x-\ox\|\\
&\leq q-r+(\vaep\|F\|+\eta)(r-q)+(\vaep+\eta)(1+\kappa\|F\|)\|x-\ox\|\\
&\leq T_{\O}^F(x)-T_{\O}^F(\ox)+[\kappa(\vaep\|F\|+\eta)+(\vaep+\eta)(1+\kappa\|F\|)]\|x-\ox\|,
\end{aligned}$$
which in turn implies that
$$\liminf_{x\to\ox}\frac{T_{\O}^F(x)-T_{\O}^F(\ox)-\la x^*,x-\ox\ra}{\|x-\ox\|}\geq -(1+2\kappa\|F\|)\vaep= -\ell\vaep,$$
since $\eta>0$ was chosen arbitrarily. Thus we get (\ref{eq1Theo2}) and complete the proof of the theorem.
$\h$

The following corollary gives us a exact characterization for Fr${\rm \acute{e}}$chet subdifferential of minimal time functions at points outside $\O$ obtained in \cite{SY}.

\begin{Corollary} Let $\ox\notin \O$ with $r:=T_{\O}^F(\ox)<\infty.$ Then we have
$$\Hat{\partial} T_{\O}^F(\ox)=\Hat{N}(\ox;\O_r)\cap\{x^*\in X^*\ |\ \sigma_F(-x^*)=1\}.$$
\end{Corollary}

We recall the one-sided limiting subdifferential  for a function $g$ defined by
\begin{equation}\label{rightlimsub}
  \partial_{\geq}g(\ox):=\Limsup_{x\xrightarrow{g^+}\ox,\vaep\downarrow 0}\Hat{\partial}_{\vaep}g(\ox),
\end{equation}
where the symbol $x\xrightarrow{g^+}\ox$ signifies that $x\to \ox$ with $g(x)\to g(\ox)$ and $g(x)\geq g(\ox).$

We also recall that a function $\ph: X^*\to\overline{\mathbb R}$ is sequentially ${\rm weak^*}$ continuous at $x^*$ if for any
sequence $x_k^*\xrightarrow{w^*}x^*$ we have $\ph(x_k^*)\to \ph(x^*)$ as $k\to\infty.$

 Removing the calmness assumption, the following theorem gives us a presentation of the limiting subdifferential  of the minimal time function at points outside $\O$, which improves the result of \cite[Theorem 6.5]{singbmn}.

\begin{Theorem}\label{Theosing1} Let $X$ be a Banach space and  let $\ox\notin\O$ with $T_{\O}^F(\ox)=r<\infty.$  Suppose that $T_{\O}^F$ is continuous around $\ox,$ and  function $\sigma_F$ is sequentially $weak^*$ continuous at every point on the set $-[N(\ox;\O_r)]\cap S^*.$ Then we have
\begin{equation}\label{eq1Theosing1}
 \partial_{\geq} T_{\O}^F(\ox)= N(\ox;\O_r)\cap S^*.
\end{equation}
\end{Theorem}
{\bf Proof.}
The inclusion $``\subset"$ in (\ref{eq1Theosing1}) follows from \cite[Theorem 6.5]{singbmn}, which requires that $X$ is a Banach space. To justify the opposite inclusion $``\supset "$ therein, fix any $x^*\in N(\ox;\O_r)\cap S^*$ and find sequences $\vaep_k\downarrow0,x_k\xrightarrow{\O_r}\ox,$ and $x_k^*\xrightarrow{w^*}x^*$ as $k\to\infty$ with $x_k^*\in\Hat N_{\vaep_k}(x_k;\O_r),$ $k\in\mathbb N.$ The sequential $\text{weak}^*$ continuity of $\sigma_F$ at $-x^*$ ensures that
$$\gamma_k:=\sigma_F(-x^*_k)\to \sigma_F(-x^*)=1\ \text{as}\ k\to\infty.$$
By the definition of $S^*$ we may assume without loss of generality that
\begin{equation}\label{eq2Theosing1} \frac{x_k^*}{\gamma_k}\in \Hat N_{\vaep_k/\gamma_k}(x_k;\O_r)\cap S^*\ \text{for all}\ k\in\mathbb N.
\end{equation}
It follows further that $T_{\O}^F(x_k)\geq r$ for large $k,$ since the opposite assumption on $T_{\O}^F(x_k)<r$ implies by the continuity of $T_{\O}^F$ that $x_k\in{\rm int} (\O_r),$ which contradicts the condition $x^*\ne 0.$  Employing Theorem \ref{Theoepsub}, find a sequence $\vaep_k'\downarrow0$ such that
$$\frac{x_k^*}{\gamma_k}\in \Hat{\partial}_{\vaep_k'}T_{\O}^F(x_k)\ \text{for all}\ k\in\mathbb N.$$
Passing there to the limit $k\to\infty$ justifies equality (\ref{eq1Theosing1}). $\h$


\end{document}